\documentclass{article}
\usepackage{amssymb}

\input{tcilatex}
\begin{document}

\begin{center}
\textbf{Numerical solution of one-dimensional Sine--Gordon equation using
Reproducing Kernel Hilbert Space Method}

Ali Akg\"{u}l and Mustafa Inc

Department of Mathematics, Education Faculty, Dicle University, 21280
Diyarbak\i r / TURKEY

Department of Mathematics, Science Faculty, F\i rat University, 23119 Elaz\i 
\u{g} / TURKEY

aliakgul00727@gmail.com, minc@firat.edu.tr
\end{center}

\textbf{Abstract:} In this paper, we propose a reproducing kernel Hilbert
space method (RKHSM) for solving the sine--Gordon (SG) equation with initial
and boundary conditions based on the reproducing kernel theory. Its exact
solution is represented in the form of series in the reproducing kernel
Hilbert space. Some numerical examples have been studied to demonstrate the
accuracy of the present method. The results obtained from the method are
compared with the exact solutions and the earlier works. Results of
numerical examples show that the presented method is simple and effective.


\bigskip

\textbf{1. Introduction}

Since Russell's [1-2] first observation for some special waves with
characteristic properties on a canal in 1834 many scientists among whom
Korteweg and de Vries [3], who derived the equation concerning the
propagation of waves in one direction on the free surface of a shallow canal
also known as the KdV equation, have investigated them more extensively. The
name soliton was first used in Zabusky and Kruskal [4] in order to emphasize
that a soliton is a localized entity which keeps its identity after
interaction. Equations which also lead to solitary waves are the
sine-Gordon, the cubic Schr\"{o}dinger equation, etc. Because of the great
importance of the study of physical phenomena theoretical solutions of
nonlinear and especially of soliton type equations have been developed
during the last years [5].

\bigskip

The nonlinear one-dimensional Sine--Gordon (SG) equation is a very important
nonlinear hyperbolic partial differential equation PDE. It appears in
differential geometry gained its significance because of the collisional
behaviors of solitons that arise from these equations. It was basically
considered in the nineteenth century in the course of surface of constant
negative curvature. This equation attracted a lot of attention in the 1970s
due to the presence of soliton solutions [6-7]. The sine-Gordon equation
comes out in a number of physical applications [8--10] including
applications in the chain of coupled pendulums and modelling the propagation
of transverse electromagnetic (TEM) wave on a superconductor transmission
system.

\bigskip

Consider the one-dimensional nonlinear sine-Gordon equation

\begin{equation}
\frac{\partial ^{2}u}{\partial t^{2}}(x,t)=\frac{\partial ^{2}u}{\partial
x^{2}}(x,t)-\sin (u(x,t)),\text{ \ }a\leq x\leq b,\text{ }t\geq 0,  \tag{1.1}
\end{equation}%
with initial conditions%
\begin{eqnarray}
u(x,0) &=&f(x),\text{ }a\leq x\leq b,  \TCItag{1.2} \\
\frac{\partial u}{\partial t}(x,0) &=&g(x),\text{ }a\leq x\leq b,  \nonumber
\end{eqnarray}%
and the boundary conditions%
\begin{equation}
u(a,t)=h_{1}(t),\text{ }u(b,t)=h_{2}(t),\text{ }t\geq 0.  \tag{1.3}
\end{equation}

Since the exact solution of the SG equation can only be obtained in special
situations, many numerical schemes are constructed to solve the SG equation.
The study of numerical solutions of the sine-Gordon equation have been
investigated considerably in the last few years. For solving $\left(
1.1\right) $, for instance, high-order solution of one-dimensional
sine-Gordon equation using compact finite difference and DIRKN methods [8] .
The authors of [11] proposed a numerical method form solving $\left(
1.1\right) $ using collocation and radial basis functions. Also, the
boundary integral equation approach is used in [12] . Bratsos has proposed a
numerical scheme for solving one dimensional sine-Gordon equation and a
third-order numerical scheme for the two-dimensional sine-Gordon equation in
[13-14] respectively. Also, in [15] , a numerical method using radial basis
function for the solution of two-dimensional sine-Gordon equation is used.
In addition, several authors recommended spectral methods and Fourier
pseudospectral method for solving nonlinear wave equation using a discrete
Fourier series and Chebyshev orthogonal polynomials [16--18]. Ma and Wu [19]
presented a meshless scheme by using a multiquadric (MQ) quasi-interpolation
named $L_{\varepsilon }$ without solving a large-scale linear system of
equations, but a polynomial $p(x)$ was needed to improved the accuracy of
the scheme. More recently, another numerical work has been investigated by
Khaliq et al. [20].

\bigskip

One difficult point of numerically solving the sine-Gordon equation is how
to solve the nonlinear system resulting by discretization. The Newton
iteration method is often used as a basis for designing numerical schemes.
In this case, the Jacobian matrices have to be established, inverted and
possibly updated during the iteration. When high approximation accuracy is
desired, it requires one to use a sufficiently small time-step and
sufficiently fine grids. Thus it demands a large amount of computational
effort [21].

\bigskip

In this paper, we solve Eqs. (1.1) and (1.3) by using Reproducing Kernel
Method (or RKM). The nonlinear problem is solved easily and elegantly by
using RKM. The technique has many advantages over the classical techniques.
It also avoids discretization and provides an efficient numerical solution
with high accuracy, minimal calculation, avoidance of physically unrealistic
assumptions. In the next section, we will describe the procedure.

\bigskip

The theory of reproducing kernels [22] was used for the first time at the
beginning of the 20 th century by Zaremba in his work on boundary value
problems for harmonic and biharmonic functions. Reproducing kernel theory
has important applications in numerical analysis, differential equations,
probability and statistics and so on [23-35]. Recently, using the RKM, some
authors discussed fractional differential equation, nonlinear oscillator
with discontinuity, singular nonlinear two-point periodic boundary value
problems, integral equations and nonlinear partial differential equations
and so on [23-35].

\bigskip

The paper is organized as follows. Section 2 is devoted to several
reproducing kernel spaces. S{\normalsize olution represantation in}\textbf{\ 
}$W\left( \Omega \right) $ and a linear operator have been presented in
Section 3. Section 4 provides the main results, the exact and approximate
solution of Eqs.(1.1) and (1.3) and an iterative method are developed for
the kind of problems in the reproducing kernel space. We have proved that
the approximate solution converges to the exact solution uniformly. Some
numerical experiments are illustrated in Section 5. We provide some
conclusions in the last section.

\bigskip

\textbf{2. Preliminaries}\bigskip

\textbf{2.1. Reproducing Kernel Spaces}\ \bigskip

In this section, we define some useful reproducing kernel spaces.\bigskip

\textbf{Definition 2.1.} \textit{(Reproducing kernel)}. Let $E$ be a
nonempty abstract set. A function $K:E\times E\longrightarrow C$ is a
reproducing kernel of the Hilbert space $H$ if and only if

a) $\forall t\in E,$ $K\left( .,t\right) \in H,$

b) $\forall t\in E,$ $\forall \varphi \in H,$ $\left\langle \varphi \left(
.\right) ,K\left( .,t\right) \right\rangle =\varphi \left( t\right) .$%
\bigskip

The last condition is called "the reproducing property"\ the value of the
function $\varphi $ at the point $t$ is reproduced by the inner product of $%
\varphi $ with $K\left( .,t\right) $.\bigskip

\bigskip

\textbf{Definition 2.2.}

\[
W_{2}^{3}[0,1]=\left\{ 
\begin{array}{c}
u(x)\mid u(x),\text{ }u^{\prime }(x),\text{ }u^{\prime \prime }(x)\text{\ \
are absolutely continuous real value functions in }[0,1], \\ 
\\ 
u^{\left( 3\right) }(x)\in L^{2}[0,1],\text{ }x\in \lbrack 0,1],\text{ }%
u(0)=0,\text{ }u(1)=0,\text{ }%
\end{array}%
\right\} , 
\]%
The inner product and the norm in $W_{2}^{3}[0,1]$ are defined respectively
by%
\[
\left\langle u(x),g(x)\right\rangle _{{\LARGE W}_{{\LARGE 2}}^{{\LARGE 3}%
}}=u(0)g(0)+u^{\prime }(0)g^{\prime }(0)+u^{\prime }(1)g^{\prime
}(1)+\int_{0}^{4}u^{(3)}(x)g^{(3)}(x)dx,\text{ \ }u(x),\text{ }g(x)\in
W_{2}^{3}[0,1], 
\]%
and%
\[
\left\Vert u\right\Vert _{{\LARGE W}_{{\LARGE 2}}^{{\LARGE 3}}}=\sqrt{%
\left\langle u,u\right\rangle _{_{{\LARGE W}_{{\LARGE 2}}^{{\LARGE 3}}}}},\
u\in W_{2}^{3}[0,1]. 
\]%
The space $W_{2}^{3}[0,1]$ \ is a reproducing kernel space, i.e., for each
fixed $y\in \lbrack 0,1]$ \ and any $u(x)\in W_{2}^{3}[0,1],$ there exists a
function $R_{y}(x)$ such that%
\[
u(y)=\left\langle u(x),\ R_{y}(x)\right\rangle _{W_{2}^{3}}. 
\]

\textbf{Definition 2.3.} \ 
\[
W_{2}^{3}[0,T]=\left\{ 
\begin{array}{c}
v(t)\mid u(t),\text{ }u^{\prime }(t),u^{\prime \prime }(t)\text{ are
absolutely continuous in \ }[0,T]\text{,} \\ 
\\ 
u^{(3)}(t)\in L^{2}[0,1],\text{ }t\in \lbrack 0,T],\text{ }u(0)=0,\text{ }%
u^{\prime }(0)=0.%
\end{array}%
\right\} , 
\]%
The inner product and the norm in $W_{2}^{3}[0,T]$ \ are defined
respectively by%
\[
\left\langle u(t),g(t)\right\rangle
_{_{W_{2}^{3}}}=\sum_{i=0}^{2}u^{(i)}(0)g^{(i)}(0)+%
\int_{0}^{1}u^{(3)}(t)g^{(3)}(t)dt,\text{ \ }u(t),g(t)\in W_{2}^{3}[0,T], 
\]%
and%
\[
\left\Vert u\right\Vert _{W_{2}^{3}}=\sqrt{\left\langle u,u\right\rangle
_{_{W_{2}^{3}}}},\text{ \ \ }u\in W_{2}^{3}[0,T]. 
\]%
The space $W_{2}^{3}[0,T]$ \ is a reproducing kernel space and its
reproducing kernel function ${\normalsize r}_{s}{\normalsize (t)}$ is given
by%
\[
r_{s}\left( t\right) =\left\{ 
\begin{array}{c}
\frac{1}{4}s^{2}t^{2}+\frac{1}{12}s^{2}t^{3}-\frac{1}{24}st^{4}+\frac{1}{120}%
t^{5},\text{ \ \ }t\leq s, \\ 
\\ 
\frac{1}{4}s^{2}t^{2}+\frac{1}{12}s^{3}t^{2}-\frac{1}{24}ts^{4}+\frac{1}{120}%
s^{5},\text{ \ \ }t>s.%
\end{array}%
\right. 
\]

\bigskip

\textbf{Definition 2.4.}\ 

\[
W_{2}^{1}[0,1]=\left\{ 
\begin{array}{c}
u(x)\mid u(x)\text{ is absolutely continuous in \ }[0,1]\text{,} \\ 
u^{\prime }(x)\in L^{2}[0,1],\text{ }x\in \lbrack 0,1].%
\end{array}%
\right\} , 
\]%
The inner product and the norm in $W_{2}^{1}[0,1]$ \ are defined
respectively by%
\[
\left\langle u(x),g(x)\right\rangle _{{\LARGE W}_{{\LARGE 2}}^{{\LARGE 1}%
}}=u(0)g(0)+\int_{0}^{1}u^{\prime }(x)g^{\prime }(x)dx,\text{ \ }%
u(x),g(x)\in W_{2}^{1}[0,1], 
\]%
and%
\[
\left\Vert u\right\Vert _{{\LARGE W}_{{\LARGE 2}}^{{\LARGE 1}}}=\sqrt{%
\left\langle u,u\right\rangle _{_{{\LARGE W}_{{\LARGE 2}}^{{\LARGE 1}}}}},%
\text{ \ \ }u\in W_{2}^{1}[0,1]. 
\]%
The space $W_{2}^{1}[0,1]$ \ is a reproducing kernel space and its
reproducing kernel function ${\normalsize Q}_{y}{\normalsize (x)}$ is given
by%
\[
{\normalsize Q}_{y}{\normalsize (x)}{\LARGE =}\left\{ 
\begin{array}{c}
1+x,\text{ \ \ }x\leq y, \\ 
\\ 
1+y,\text{ \ \ }x>y.%
\end{array}%
\right. 
\]

\bigskip

\textbf{Definition 2.5.}

\[
W_{2}^{1}[0,T]=\left\{ 
\begin{array}{c}
u(t)\mid u(t)\text{\ is absolutely continuous in \ }[0,T]\text{,} \\ 
\\ 
u^{\prime }(t)\in L^{2}[0,T],\text{ }t\in \lbrack 0,T].%
\end{array}%
\right\} , 
\]%
The inner product and the norm in $W_{2}^{1}[0,T]$\ are defined respectively
by%
\[
\left\langle u(t),g(t)\right\rangle _{_{{\LARGE W}_{{\LARGE 2}}^{{\LARGE 1}%
}}}=u(0)g(0)+\int_{0}^{T}u^{\prime }(t)g^{\prime }(t)dt,\text{ \ }%
u(t),g(t)\in W_{2}^{1}[0,T], 
\]%
and%
\[
\left\Vert u\right\Vert _{{\LARGE W}_{{\LARGE 2}}^{{\LARGE 1}}}=\sqrt{%
\left\langle u,u\right\rangle _{{\LARGE W}_{{\LARGE 2}}^{{\LARGE 1}}}},\text{
\ \ }u\in W_{2}^{1}[0,T]. 
\]%
The space $W_{2}^{1}[0,T]$\ \ is a reproducing kernel space and its
reproducing kernel function $q_{s}(t)$ is given by%
\[
q_{s}(t)=\left\{ 
\begin{array}{c}
1+t,\text{ \ \ }t\leq s, \\ 
\\ 
1+s,\text{ \ \ }t>s.%
\end{array}%
\right. 
\]

\bigskip

\textbf{Theorem 2.1. }The space $W_{2}^{3}\left[ 0,1\right] $ is a complete
reproducing kernel space whose reproducing kernel function $R_{y}\left(
x\right) $ is given as,%
\begin{equation}
R_{y}\left( x\right) =\left\{ 
\begin{array}{c}
\sum_{i=1}^{6}c_{i}\left( y\right) x^{i-1},\quad x\leq y, \\ 
\\ 
\sum_{i=1}^{6}d_{i}\left( y\right) x^{i-1},\quad x>y.%
\end{array}%
\right.   \tag{2.1}
\end{equation}%
where%
\[
c_{1}(y)=0,
\]%
\[
{\normalsize c}_{2}{\normalsize (y)=-}\frac{1}{122}y^{5}+\frac{5}{244}y^{4}-%
\frac{127}{244}y^{2}+\frac{31}{61}y,
\]%
\[
{\normalsize c}_{3}{\normalsize (y)=-}\frac{1}{2928}y^{5}+\frac{127}{5856}%
y^{4}-\frac{1}{12}y^{3}+\frac{1137}{1952}y^{2}-\frac{127}{244}y,
\]%
\[
{\normalsize c}_{4}{\normalsize (y)=0},
\]%
\[
{\normalsize c}_{5}{\normalsize (y)=}\frac{1}{2938}y^{5}-\frac{5}{5856}y^{4}+%
\frac{127}{5856}y^{2}-\frac{31}{1464}y,
\]%
\[
{\normalsize c}_{6}{\normalsize (y)=-}\frac{1}{7320}y^{5}+\frac{1}{2928}%
y^{4}-\frac{1}{2928}y^{2}-\frac{1}{122}y+\frac{1}{120},
\]%
\[
d_{1}(y)=\frac{1}{120}y^{5},
\]%
\[
{\normalsize d}_{2}{\normalsize (y)=-}\frac{1}{122}y^{5}-\frac{31}{1464}%
y^{4}-\frac{127}{244}y^{2}+\frac{31}{61}y,
\]%
\[
{\normalsize d}_{3}{\normalsize (y)=-}\frac{1}{2928}y^{5}+\frac{127}{5856}%
y^{4}+\frac{1137}{1952}y^{2}-\frac{127}{244}y,
\]%
\[
{\normalsize d}_{4}{\normalsize (y)=-}\frac{1}{12}y^{2},
\]%
\ 
\[
{\normalsize d}_{5}{\normalsize (y)=}\frac{1}{2928}y^{5}-\frac{5}{5856}y^{4}+%
\frac{127}{5856}y^{2}+\frac{5}{244}y,
\]%
\[
{\normalsize d}_{6}{\normalsize (y)=-}\frac{1}{7320}y^{5}+\frac{1}{2928}%
y^{4}-\frac{1}{2928}y^{2}-\frac{1}{122}y.
\]

\textbf{Proof: }Let $u\in W_{2}^{3}[0,1]$ and $0\leq y\leq 1.$ Define $R_{y}$
by (2.1). Note that%
\[
R_{y}^{\prime }\left( x\right) =\left\{ 
\begin{array}{c}
\sum_{i=1}^{5}ic_{i+1}\left( y\right) x^{i-1},\quad x<y, \\ 
\\ 
\sum_{i=1}^{5}id_{i+1}\left( y\right) x^{i-1},\quad x>y,%
\end{array}%
\right. 
\]

\[
R_{y}^{\prime \prime }\left( x\right) =\left\{ 
\begin{array}{c}
\sum_{i=1}^{4}i(i+1)c_{i+2}\left( y\right) x^{i-1},\quad x<y, \\ 
\\ 
\sum_{i=1}^{4}i(i+1)d_{i+2}\left( y\right) x^{i-1},\quad x>y,%
\end{array}%
\right. 
\]

\[
R_{y}^{(3)}\left( x\right) =\left\{ 
\begin{array}{c}
\sum_{i=1}^{3}i(i+1)(i+2)c_{i+3}\left( y\right) x^{i-1},\quad x<y, \\ 
\\ 
\sum_{i=1}^{3}i(i+1)(i+2)d_{i+3}\left( y\right) x^{i-1},\quad x>y,%
\end{array}%
\right. 
\]

\[
R_{y}^{(4)}\left( x\right) =\left\{ 
\begin{array}{c}
\sum_{i=1}^{2}i(i+1)(i+2)(i+3)c_{i+4}\left( y\right) x^{i-1},\quad x<y, \\ 
\\ 
\sum_{i=1}^{2}i(i+1)(i+2)(i+3)d_{i+4}\left( y\right) x^{i-1},\quad x>y,%
\end{array}%
\right. 
\]

and%
\[
R_{y}^{(5)}\left( x\right) =\left\{ 
\begin{array}{c}
120c_{6}(y),\quad x<y, \\ 
\\ 
120d_{6}(y),\quad x>y.%
\end{array}%
\right. 
\]

By Definition 2.2 and integrating by parts two times, we obtain%
\begin{equation}
\left. 
\begin{array}{c}
\left\langle u(x),{\normalsize R}_{y}{\normalsize (x)}\right\rangle _{_{%
{\LARGE W}_{2}^{3}}}=u(0)R_{y}(0)+u^{\prime }(0)R_{y}^{\prime }(0)+u^{\prime
}(1)R_{y}^{\prime }(1)+{\normalsize u}^{\prime \prime }(1){\normalsize R}%
_{y}^{(3)}{\normalsize (1)} \\ 
\text{ \ \ \ \ \ \ \ \ \ \ }+{\normalsize u}(1){\normalsize R}_{y}^{(5)}%
{\normalsize (1)}-{\normalsize u}(0){\normalsize R}_{y}^{(5)}{\normalsize %
(0)-}\int_{0}^{1}u{\normalsize (x)R}_{y}^{(6)}{\normalsize (x)dx} \\ 
\\ 
\text{ \ \ \ \ \ \ \ \ }=u^{\prime }(0)(R_{y}^{\prime }(0)+{\normalsize R}%
_{y}^{(4)}{\normalsize (0)})+u^{\prime }(1)(R_{y}^{\prime }(1)-{\normalsize R%
}_{y}^{(4)}{\normalsize (1)}) \\ 
+{\normalsize u}^{\prime \prime }(1){\normalsize R}_{y}^{(3)}{\normalsize (1)%
}-{\normalsize u}^{\prime \prime }(0){\normalsize R}_{y}^{(3)}{\normalsize %
(0)} \\ 
{\normalsize +}\int_{0}^{y}{\normalsize R}_{y}^{(5)}{\normalsize (x)}u%
{\normalsize (x)dx+}\int_{y}^{1}{\normalsize R}_{y}^{(5)}{\normalsize (x)}u%
{\normalsize (x)dx} \\ 
\\ 
=u^{\prime }(0)(c_{2}(y)+{\normalsize 24c}_{5}(y))-{\normalsize u}^{\prime
\prime }(0){\normalsize (6c}_{4}(y){\normalsize )} \\ 
\text{ \ \ \ \ \ \ \ \ \ \ \ \ \ \ }+u^{\prime
}(1)(d_{2}(y)+2d_{3}(y)+3d_{4}(y)-20d_{5}(y)-115d_{6}(y)) \\ 
+{\normalsize u}^{\prime \prime }(1)(6d_{4}(y)+24d_{5}(y)+60d_{6}(y)) \\ 
{\normalsize +}\int_{0}^{y}120c_{6}(y)u{\normalsize (x)dx+}%
\int_{y}^{1}120d_{6}(y)u{\normalsize (x)dx} \\ 
\\ 
=120u(y)(\frac{1}{120})=u(y).%
\end{array}%
\right.   \nonumber
\end{equation}

\textbf{Definition 2.6.}%
\[
{\normalsize W}\left( \Omega \right) {\LARGE =}\left\{ 
\begin{array}{c}
u(x,t)\mid \frac{{\LARGE \partial }^{4}{\LARGE u}}{{\LARGE \partial x}^{2}%
{\LARGE \partial t}^{2}},\text{ is completely continuos in }\Omega
=[0,1]\times \lbrack 0,T], \\ 
\\ 
\frac{{\LARGE \partial }^{6}{\LARGE u}}{{\LARGE \partial x}^{3}{\LARGE %
\partial t}^{3}}\in L^{2}\left( \Omega \right) ,\text{ }u(x,0)=0,\text{ }%
\frac{{\LARGE \partial u(x,0)}}{{\LARGE \partial t}}=0,\text{ }u(0,t)=0,%
\text{ }u(1,t)=0.%
\end{array}%
\right\} ,
\]%
The inner product and the norm in ${\normalsize W}\left( \Omega \right) $
are defined respectively by%
\begin{eqnarray*}
\left\langle u(x,t),g(x,t)\right\rangle _{_{W}} &=&\sum_{i=0}^{2}\int_{0}^{T}
\left[ \frac{\partial ^{3}}{\partial t^{3}}\frac{\partial ^{i}}{\partial
x^{i}}u(0,t)\frac{\partial ^{3}}{\partial t^{3}}\frac{\partial ^{i}}{%
\partial x^{i}}g(0,t)\right] dt \\
&& \\
&&+\sum_{j=0}^{2}\left\langle \frac{\partial ^{j}}{\partial t^{j}}u(x,0),%
\frac{\partial ^{j}}{\partial t^{j}}g(x,0)\right\rangle _{{\LARGE W}_{%
{\LARGE 2}}^{{\LARGE 3}}} \\
&& \\
&&+\int_{0}^{1}\int_{0}^{T}\left[ \frac{\partial ^{3}}{\partial x^{3}}\frac{%
\partial ^{3}}{\partial t^{3}}u(x,t)\frac{\partial ^{3}}{\partial x^{3}}%
\frac{\partial ^{3}}{\partial t^{3}}g(x,t)\right] dtdx,
\end{eqnarray*}%
and%
\[
\left\Vert u\right\Vert _{W}=\sqrt{\left\langle u,u\right\rangle _{W}},\
u\in W\left( \Omega \right) .
\]

\bigskip

\textbf{Theorem 2.2. }$W(\Omega )$ is a reproducing kernel space and its
reproducing kernel function is%
\[
K_{\left( y,s\right) }\left( x,t\right) =R_{y}(x)r_{s}(t), 
\]%
such that for any $u(x,t)\in W\left( \Omega \right) ,$%
\[
u(y,s)=\left\langle u(x,t),K_{\left( y,s\right) }\left( x,t\right)
\right\rangle _{W}, 
\]%
and%
\[
K_{\left( y,s\right) }\left( x,t\right) =K_{\left( x,t\right) }\left(
y,s\right) . 
\]

\bigskip

\textbf{Definition 2.7.}%
\[
\widehat{W}\left( \Omega \right) {\LARGE =}\left\{ 
\begin{array}{c}
u(x,t)\mid u(x,t)\text{ is completely continuos in \ }\Omega =[0,1]\times
\lbrack 0,T]\text{,} \\ 
\frac{{\LARGE \partial }^{{\LARGE 2}}{\huge u}}{{\huge \partial x\partial t}}%
\in L^{2}\left( \Omega \right) .%
\end{array}%
\right\} {\LARGE .} 
\]%
The inner product and the norm in $\ \widehat{W}\left( \Omega \right) $ are
defined respectively by

\begin{eqnarray*}
\left\langle u(x,t),g(x,t)\right\rangle _{\widehat{W}} &=&\int_{0}^{T}\left[ 
\dfrac{\partial }{\partial t}u(0,t)\dfrac{\partial }{\partial t}g(0,t)\right]
dt \\
&& \\
&&+\left\langle u(x,0),g(x,0\right\rangle _{{\LARGE W}_{{\LARGE 2}}^{{\LARGE %
1}}} \\
&& \\
&&+\int_{0}^{1}\int_{0}^{T}\left[ \frac{\partial }{\partial x}\frac{\partial 
}{\partial t}u(x,t)\frac{\partial }{\partial x}\frac{\partial }{\partial t}%
g(x,t)\right] dtdx,
\end{eqnarray*}%
and

\[
\left\Vert u\right\Vert _{\widehat{W}}=\sqrt{\left\langle u,u\right\rangle _{%
\widehat{W}}},\text{ \ }u\in \widehat{W}\left( \Omega \right) . 
\]%
$\widehat{W}\left( \Omega \right) $\ is a reproducing kernel space and its
reproducing kernel function $G_{\left( y,s\right) }\left( x,t\right) $ is%
\[
G_{\left( y,s\right) }\left( x,t\right) =Q_{y}(x)q_{s}(t). 
\]

\bigskip

\textbf{3. Solution represantation in }$W\left( \Omega \right) .$

In this section, the solution of equation (1.1) is given in the reproducing
kernel space $W\left( \Omega \right) $. On defining the linear operator $%
L:W\left( \Omega \right) \rightarrow \widehat{W}\left( \Omega \right) $\ as

\[
Lv=\frac{\partial ^{2}v}{\partial t^{2}}(x,t)-\frac{\partial ^{2}v}{\partial
x^{2}}(x,t).
\]%
After homogenizing the initial and boundary conditions model problem (1.1)
changes the following problem:

\begin{equation}
\left\{ 
\begin{array}{c}
Lv=M(x,\text{ }t,\text{ }v(x,\text{ }t)),\text{ }\left( x,t\right) \in
\lbrack 0,1]\times \lbrack 0,T], \\ 
\\ 
v(x,0)=\frac{{\LARGE \partial v}}{{\LARGE \partial t}}(x,0)=v(0,t)=v(1,t)=0.%
\end{array}%
\right. \text{\ }  \tag{3.1}
\end{equation}%
We replace $v(x,t)$ with $u(x,t)$ in (3.1), for simplicity.

\bigskip

\textbf{Lemma 3.1. }$L$ is a bounded linear operator.

\bigskip

\textbf{Proof:}%
\begin{eqnarray*}
\left\Vert Lu\right\Vert _{\widehat{W}}^{2} &=&\int_{0}^{T}\left[ \frac{%
\partial }{\partial t}Lu(0,t)\right] ^{2}dt+\left\langle
Lu(x,0),Lu(x,0)\right\rangle _{{\huge W}_{{\LARGE 2}}^{{\LARGE 1}}} \\
&& \\
&&+\int_{0}^{1}\int_{0}^{T}\left[ \frac{\partial }{\partial x}\frac{\partial 
}{\partial t}Lu(x,t)\right] ^{2}dxdt \\
&& \\
&=&\int_{0}^{T}\left[ \frac{\partial }{\partial t}Lu(0,t)\right] ^{2}dt+%
\left[ Lu(0,0)\right] ^{2} \\
&& \\
&&+\int_{0}^{1}\left[ \frac{\partial }{\partial x}Lu(x,0)\right]
^{2}dx+\int_{0}^{1}\int_{0}^{T}\left[ \frac{\partial }{\partial x}\frac{%
\partial }{\partial t}Lu(x,t)\right] ^{2}dxdt,
\end{eqnarray*}%
since%
\begin{eqnarray*}
u(x,t) &=&\left\langle u(\xi ,\eta ),K_{\left( {\LARGE x,t}\right) }\left(
\xi ,\eta \right) \right\rangle _{W}, \\
Lu(x,t) &=&\left\langle u(\xi ,\eta ),LK_{\left( x,t\right) }\left( \xi
,\eta \right) \right\rangle _{W},
\end{eqnarray*}%
from the continuity of $K_{\left( {\LARGE x,t}\right) }\left( \xi ,\eta
\right) $, we have%
\[
\left\vert Lu(x,t)\right\vert \leq \left\Vert u\right\Vert _{W}\left\Vert
LK_{\left( x,t\right) }\left( \xi ,\eta \right) \right\Vert _{W}\leq
a_{0}\left\Vert u\right\Vert _{W}.
\]%
similarly for $i=0,1$%
\[
\frac{\partial ^{i}}{\partial x^{i}}Lu(x,t)=\left\langle u(\xi ,\eta ),\frac{%
\partial ^{i}}{\partial x^{i}}{\normalsize LK}_{\left( {\LARGE x,t}\right)
}\left( \xi ,\eta \right) \right\rangle _{W}{\LARGE ,}
\]%
\[
\frac{\partial }{\partial t}\frac{\partial ^{i}}{\partial x^{i}}%
Lu(x,t)=\left\langle u(\xi ,\eta ),\frac{\partial }{\partial t}\frac{%
\partial ^{i}}{\partial x^{i}}LK_{\left( {\LARGE x,t}\right) }\left( \xi
,\eta \right) \right\rangle _{W}{\LARGE ,}
\]%
and then%
\[
\left\vert \frac{\partial ^{i}}{\partial x^{i}}Lu(x,t)\right\vert \leq e%
\text{{\large $_{i}$}}\left\Vert u\right\Vert _{W}{\normalsize ,}
\]%
\[
\left\vert \frac{\partial }{\partial t}\frac{\partial ^{i}}{\partial x^{i}}%
Lu(x,t)\right\vert \leq f\text{{\large $_{i}$}}\left\Vert u\right\Vert _{W}%
{\normalsize .}
\]%
Therefore%
\[
\left\Vert Lu(x,t)\right\Vert _{\widehat{W}}^{2}\text{ }{\normalsize \leq }%
\sum_{i=0}^{1}\left( e_{i}^{2}+f_{i}^{2}\right) \left\Vert u\right\Vert
_{W}^{2}\leq {\normalsize a}^{2}\left\Vert u\right\Vert _{W}^{2}{\LARGE .}%
\text{ \ \ \ \ }\square 
\]

Now, choose\ a countable dense subset \ \ $\left\{ \left( x_{1},t_{1}\right)
,\left( x_{2},t_{2}\right) ,...\right\} $ in $\Omega =[0,1]\times \lbrack
0,T]$ \ and define%
\[
\varphi _{i}(x,t)=G_{\left( x_{i},t_{i}\right) }\left( x,t\right) ,\text{ \ }%
\Psi _{i}(x,t)=L^{\ast }\varphi _{i}(x,t), 
\]%
where $L^{\ast }$ is the adjoint operator of $L.$ The orthonormal system $%
\left\{ \widehat{\Psi }_{i}(x,t)\right\} _{i=1}^{\infty }$ of $W\left(
\Omega \right) $ \ can be derived from the process of Gram-Schmidt
orthogonalization of $\left\{ \Psi _{i}(x,t)\right\} _{i=1}^{\infty }$ as%
\[
\widehat{\Psi }_{i}(x,t)=\sum_{k=1}^{i}\beta _{ik}\Psi _{k}(x,t). 
\]

\bigskip

\textbf{Theorem 3.1. }Suppose that $\left\{ (x_{i},t_{i})\right\}
_{i=1}^{\infty }$ is dense in $\Omega ;$ then $\left\{ \Psi
_{i}(x,t)\right\} _{i=1}^{\infty }$ is complete system in $W\left( \Omega
\right) $ and%
\[
\Psi _{i}(x,t)=\left. L_{\left( {\LARGE y,s}\right) }K_{\left( {\LARGE y,s}%
\right) }\left( x,t\right) \right\vert _{\left( {\LARGE y,s}\right) =\left( 
{\LARGE x}_{i}{\LARGE ,t}_{i}\right) }. 
\]

\textbf{Proof: }We have%
\begin{eqnarray*}
\Psi _{i}(x,t) &=&\left( L^{\ast }\varphi _{i}\right) \left( x,t\right)
=\left\langle \left( L^{\ast }\varphi _{i}\right) \left( y,s\right)
,K_{\left( {\LARGE x,}\text{ }{\LARGE t}\right) }\left( y,s\right)
\right\rangle _{W} \\
&=&\left\langle \varphi _{i}\left( y,s\right) ,\text{ }L_{\left( {\LARGE y,s}%
\right) }K_{\left( {\LARGE x,t}\right) }\left( y,s\right) \right\rangle _{%
\widehat{W}} \\
&=&\left. L_{\left( {\LARGE y,s}\right) }K_{\left( {\LARGE x,}\text{ }%
{\LARGE t}\right) }\left( y,s\right) \right\vert _{\left( {\LARGE y,s}%
\right) =\left( {\LARGE x}_{{\LARGE i}},\text{ }{\LARGE t}_{{\LARGE i}%
}\right) } \\
&=&\left. L_{\left( {\LARGE y,s}\right) }K_{\left( {\LARGE y,s}\right)
}\left( x,t\right) \right\vert _{\left( {\LARGE y,}\text{ }{\LARGE s}\right)
=\left( {\LARGE x}_{{\LARGE i}}{\LARGE ,}\text{ \ }{\LARGE t}_{{\LARGE i}%
}\right) }.
\end{eqnarray*}%
Clearly $\Psi _{i}(x,t)\in W\left( \Omega \right) .$ For each fixed $%
u(x,t)\in W\left( \Omega \right) ,$ if%
\[
\left\langle u(x,t),\Psi _{i}(x,t)\right\rangle _{W}=0,\text{ \ \ }i=1,2,... 
\]%
then%
\begin{eqnarray*}
\left\langle u(x,t),\left( L^{\ast }\varphi _{i}\right) \left( x,t\right)
\right\rangle _{W} &=&\left\langle Lu(x,t),\varphi _{i}\left( x,t\right)
\right\rangle _{\widehat{W}} \\
&=&\left( {\LARGE Lu}\right) \left( {\LARGE x}_{{\LARGE i}}{\LARGE ,}\text{ }%
{\LARGE t}_{{\LARGE i}}\right) =0,\text{ \ }i=1,2,...\mathbf{.}
\end{eqnarray*}%
Note that$\ \left\{ (x_{i},t_{i})\right\} _{i=1}^{\infty }$ is dense in $%
\Omega $, hence, $\left( Lu\right) \left( x,t\right) =0.$ It follows that $%
u=0$ from the existence of $L^{-1}$. So the proof is complete. \ \ \ \ $%
\square $

\bigskip

\textbf{Theorem 3.2.} If $\left\{ (x_{i},t_{i})\right\} _{i=1}^{\infty }$ is
dense in $\Omega $, then the solution of (3.1) is

\begin{equation}
u(x,t)=\sum_{i=1}^{\infty }\sum_{k=1}^{i}\beta _{ik}M(x_{k},\text{ }t_{k},%
\text{ }u(x_{k},\text{ }t_{k}))\widehat{\Psi }_{i}(x,t).  \tag{3.2}
\end{equation}

\bigskip

\textbf{Proof. }Since\textbf{\ }$\left\{ \Psi _{i}(x,t)\right\}
_{i=1}^{\infty }$ is complete system in $W\left( \Omega \right) ,$ we have

\begin{eqnarray*}
u(x,t) &=&\sum_{i=1}^{\infty }\left\langle u(x,t),\widehat{\Psi }%
_{i}(x,t)\right\rangle _{W}\widehat{\Psi }_{i}(x,t) \\
&& \\
&=&\sum_{i=1}^{\infty }\sum_{k=1}^{i}\beta _{ik}\left\langle u(x,t),\Psi
_{k}(x,t)\right\rangle _{W}\widehat{\Psi }_{i}(x,t) \\
&& \\
&=&\sum_{i=1}^{\infty }\sum_{k=1}^{i}\beta _{ik}\left\langle u(x,t),L^{\ast
}\varphi _{k}(x,t)\right\rangle _{W}\widehat{\Psi }_{i}(x,t)
\end{eqnarray*}

\begin{eqnarray*}
\text{ \ \ \ \ \ \ \ \ \ \ \ \ \ } &=&\sum_{i=1}^{\infty
}\sum_{k=1}^{i}\beta _{ik}\left\langle Lu(x,t),\varphi
_{k}(x,t)\right\rangle _{\widehat{W}}\widehat{\Psi }_{i}(x,t) \\
&& \\
&=&\sum_{i=1}^{\infty }\sum_{k=1}^{i}\beta _{ik}\left\langle
Lu(x,t),G_{\left( {\LARGE x}_{{\LARGE k}}{\LARGE ,}\text{ }{\LARGE t}_{%
{\LARGE k}}\right) }(x,t)\right\rangle _{\widehat{W}}\widehat{\Psi }_{i}(x,t)
\\
&& \\
&=&\sum_{i=1}^{\infty }\sum_{k=1}^{i}\beta _{ik}Lu\left( x_{k},t_{k}\right) 
\widehat{\Psi }_{i}(x,t) \\
&& \\
&=&\sum_{i=1}^{\infty }\sum_{k=1}^{i}\beta _{ik}M(x_{k},\text{ }t_{k},\text{ 
}u(x_{k},\text{ }t_{k}))\widehat{\Psi }_{i}(x,t).
\end{eqnarray*}

\bigskip

Now the approximate solution $u_{n}(x,t)$ \ can be obtained from the $n-$%
term intercept of the exact solution $u(x,t)$ \ and%
\[
{\normalsize u}_{n}{\normalsize (x,t)=}\sum_{i=1}^{n}\sum_{k=1}^{i}\beta
_{ik}M(x_{k},\text{ }t_{k},\text{ }u(x_{k},\text{ }t_{k}))\widehat{\Psi }%
_{i}(x,t). 
\]%
Obviously%
\[
\left\Vert u_{n}(x,t)-u(x,t)\right\Vert \rightarrow 0,\ \ \left(
n\rightarrow \infty \right) .\text{ \ \ \ }\square 
\]

\bigskip

\textbf{4. The method implementation}\bigskip

If we write%
\[
A_{i}=\sum_{k=1}^{i}\beta _{ik}M(x_{k},\text{ }t_{k},\text{ }u(x_{k},\text{ }%
t_{k})), 
\]%
then (3.2) can be written as%
\[
u(x,t)=\sum_{i=1}^{\infty }A_{i}\widehat{\Psi }_{i}(x,t). 
\]

Now let $(x_{1},t_{1})=0;$ then from the initial and boundary conditions of
(3.1), $u(x_{1},t_{1})$ is known. We put $u_{0}(x_{1},t_{1})=u(x_{1},t_{1})$
\ and define the $n-$ term approximation to $u(x,t)$ by%
\begin{equation}
u_{n}(x,t)=\sum_{i=1}^{n}B_{i}\widehat{\Psi }_{i}(x,t).  \tag{3.3}
\end{equation}%
where%
\begin{equation}
B_{i}=\sum_{k=1}^{i}\beta _{ik}M(x_{k},\text{ }t_{k},\text{ }u_{k-1}(x_{k},%
\text{ }t_{k})).  \tag{3.4}
\end{equation}

In the sequel, we verify that the approximate solution ${\normalsize u}_{n}%
{\normalsize (x,t)}$ converges to the exact solution, uniformly.

\bigskip

\textbf{Theorem 4.1. }Suppose that $\left\Vert u_{n}\right\Vert $ is a
bounded in (3.3) and (3.1) has a uniqe solution. If $\ \left\{
(x_{i},t_{i})\right\} _{i=1}^{\infty }$ $\ $is dense in $W\left( \Omega
\right) ,$ ~then the $n-$ term approximate solution $\ u_{n}(x,t)$ can be\
derived from the above method converges to the analytical solution $u(x,t)$
of (3.1) and%
\[
u(x,t)=\sum_{i=1}^{\infty }B_{i}\widehat{\Psi }_{i}(x,t), 
\]%
where $B_{i}$ is given by (3.4).\bigskip

\textbf{Proof: }First, we prove the convergence of $u_{n}(x,t).$ From (3.3),
we infer that%
\[
{\normalsize u}_{n+1}{\normalsize (x,t)=u}_{n}{\normalsize (x,t)+B}_{n+1}%
\widehat{\Psi }_{n+1}{\normalsize (x,t),}
\]%
The orthonormality of $\ \left\{ \widehat{\Psi }_{i}\right\} _{i=1}^{\infty }
$ yields that%
\begin{equation}
\left\Vert u_{n+1}\right\Vert ^{2}{\normalsize =}\left\Vert u_{n}\right\Vert
^{2}{\normalsize +B}_{n+1}^{2}{\normalsize =}\sum_{i=1}^{n+1}{\normalsize B}%
_{i}^{2}.  \tag{3.5}
\end{equation}%
In terms of (3.5), it holds that $\left\Vert u_{n+1}\right\Vert >\left\Vert
u_{n}\right\Vert .$ Due to the condition that $\left\Vert u_{n}\right\Vert $
is bounded, $\left\Vert u_{n}\right\Vert $ is convergent and there exists a
constant $c$ such that%
\[
\sum_{i=1}^{\infty }B_{i}^{2}=c.
\]%
This implies that%
\[
\left\{ B_{i}\right\} _{i=1}^{\infty }\in l^{2}.
\]%
If $\ m>n,$ then%
\begin{eqnarray*}
\left\Vert u_{m}-u_{n}\right\Vert ^{2} &=&\left\Vert
u_{m}-u_{m-1}+u_{m-1}-u_{m-2}+...+u_{n+1}-u_{n}\right\Vert ^{2} \\
&\leq &\left\Vert u_{m}-u_{m-1}\right\Vert ^{2}+\left\Vert
u_{m-1}-u_{m-2}\right\Vert ^{2}+...+\left\Vert u_{n+1}-u_{n}\right\Vert ^{2}.
\end{eqnarray*}%
On account of%
\[
\left\Vert u_{m}-u_{m-1}\right\Vert ^{2}=B_{m}^{2},
\]%
consequently,%
\[
\left\Vert u_{m}-u_{n}\right\Vert ^{2}=\sum_{l=n+1}^{m}B_{l}^{2}\rightarrow
0,\text{ as }n\rightarrow \infty .
\]%
The completeness of $W\left( \Omega \right) $ shows that $u_{n}\rightarrow 
\widehat{u}$ as $n\rightarrow \infty $ . Now, let we prove that\ $\widehat{u}
$ \ is the solution of (3.1). Taking limits in (3.3) we get%
\[
\widehat{u}(x,t)=\sum_{i=1}^{\infty }B_{i}\widehat{\Psi }_{i}(x,t).
\]%
Note that%
\[
(L\widehat{u})(x,t)=\sum_{i=1}^{\infty }B_{i}L\widehat{\Psi }_{i}(x,t),
\]%
and%
\begin{eqnarray*}
(L\widehat{u})(x_{l},t_{l}) &=&\sum_{i=1}^{\infty }B_{i}L\widehat{\Psi }%
_{i}(x_{l},t_{l}) \\
&=&\sum_{i=1}^{\infty }B_{i}\left\langle L\widehat{\Psi }_{i}(x,t),\text{ }%
\varphi _{l}(x,t)\right\rangle _{\widehat{W}} \\
&=&\sum_{i=1}^{\infty }B_{i}\left\langle \widehat{\Psi }_{i}(x,t),\text{ }%
L^{\ast }\varphi _{l}(x,t)\right\rangle _{W} \\
&=&\sum_{i=1}^{\infty }B_{i}\left\langle \widehat{\Psi }_{i}(x,t),\text{ }%
\widehat{\Psi }_{l}(x,t)\right\rangle _{W}.
\end{eqnarray*}%
Therefore%
\begin{eqnarray*}
\sum_{l=1}^{i}\beta _{il}(L\widehat{u})(x_{l},t_{l}) &=&\sum_{i=1}^{\infty
}B_{i}\left\langle \widehat{\Psi }_{i}(x,t),\text{ }\sum_{l=1}^{i}\beta
_{il}\Psi _{l}(x,t)\right\rangle _{W} \\
&=&\sum_{i=1}^{\infty }B_{i}\left\langle \widehat{\Psi }_{i}(x,t),\text{ }%
\widehat{\Psi }_{l}(x,t)\right\rangle _{W}=B_{l}.
\end{eqnarray*}%
In wiev of (3.4), we have%
\[
{\normalsize L}\widehat{u}{\normalsize (x}_{l}{\normalsize ,t}_{l}%
{\normalsize )=M(x}_{l},\text{ }t_{l},\text{ }u({\normalsize x}_{l},\text{ }%
t_{l}))
\]%
Since $\left\{ (x_{i},t_{i})\right\} _{i=1}^{\infty }$ is dense in$\ \Omega ,
$ for each $(y,s)\in \Omega ,$ there exists a subsequence $\left\{ \left(
x_{n_{j}},t_{n_{j}}\right) \right\} _{j=1}^{\infty }$ \ such that

\[
\left( x_{n_{j}},t_{n_{j}}\right) \rightarrow (y,s),\text{ \ }\
(j\rightarrow \infty ). 
\]%
We know that%
\[
L\widehat{u}\left( x_{n_{j}},t_{n_{j}}\right) =M(x_{n_{j}},\text{ }t_{n_{j}},%
\text{ }u(x_{n_{j}},\text{ }t_{n_{j}})). 
\]%
Let $j\rightarrow \infty ;$ by the continuity of \ $f$, \ we have%
\[
(L\widehat{u})(y,s)=M(y,s,u(y,s)). 
\]%
which indicates that \ $\widehat{u}(x,t)$ satisfies (3.1).\bigskip

\textbf{Remark 4.1. }In a same manner, it can be proved that%
\[
\left\Vert \frac{\partial u_{n}(x,t)}{\partial x}-\frac{\partial u(x,t)}{%
\partial x}\right\Vert \rightarrow 0,\text{ \ as }n\rightarrow \infty , 
\]%
where%
\[
\frac{\partial u(x,t)}{\partial x}=\sum_{i=1}^{\infty }B_{i}\frac{\partial 
\widehat{\Psi }_{i}(x,t)}{\partial x}, 
\]%
and%
\[
\frac{\partial u_{n}(x,t)}{\partial x}=\sum_{i=1}^{n}B_{i}\frac{\partial 
\widehat{\Psi }_{i}(x,t)}{\partial x}, 
\]%
$B_{i}$ is given by (3.4).\bigskip

\bigskip

\textbf{5. Numerical Results}\bigskip

\textbf{Example 5.1.} In this example, we consider SG equation (1.1) without
nonlinear term $sin(u)$ in the region $\left[ 0,1\right] \times \lbrack 0,T]$%
. The initial conditions are given by%
\[
u(x,0)=\sin (\pi x),\text{ }\frac{\partial u}{\partial t}(x,0)=0, 
\]%
with the boundary conditions

\[
u(0,t)=u(1,t)=0. 
\]%
The exact solution is given in [36] as%
\[
u\left( x,t\right) =\frac{1}{2}(\sin (\pi (x+t)+\sin \pi (x-t)). 
\]

After homogenizing the initial and boundary conditions we obtain (5.1) as

\begin{equation}
\left\{ 
\begin{array}{c}
\frac{{\LARGE \partial }^{2}{\LARGE u}}{{\LARGE \partial t}^{2}}(x,t)-\frac{%
{\LARGE \partial }^{2}{\LARGE u}}{{\LARGE \partial x}^{2}}(x,t)=-\pi
^{2}\sin (\pi x), \\ 
\\ 
u(x,0)=\frac{{\LARGE \partial u}}{{\LARGE \partial t}}(x,0)=u(0,t)=u(1,t)=0.%
\end{array}%
\right.  \tag{5.1}
\end{equation}

\bigskip

\textbf{Example 5.2.} In this example, we take notice of SG equation (1.1)
in the region $\left[ 0,1\right] \times \lbrack 0,T]$. The initial
conditions are given by%
\[
u(x,0)=0,\text{ }\frac{\partial u}{\partial t}(x,0)=4\sec h(x), 
\]%
The exact solution is 
\[
u\left( x,t\right) =\arctan (\sec h(x)t). 
\]

After homogenizing the initial conditions we obtain (5.2) as

\begin{equation}
\left\{ 
\begin{array}{c}
\frac{{\LARGE \partial }^{2}{\LARGE u}}{{\LARGE \partial t}^{2}}(x,t)-\frac{%
{\LARGE \partial }^{2}{\LARGE u}}{{\LARGE \partial x}^{2}}(x,t)=-\sin
(u(x,t)+4t\sec h(x))-\frac{{\LARGE 4t}}{\cosh {\LARGE (x)}}+\frac{{\LARGE 8t}%
\sinh ^{2}{\LARGE (x)}}{\cosh ^{3}(x)}, \\ 
\\ 
u(x,0)=\frac{{\LARGE \partial u}}{{\LARGE \partial t}}(x,0)=u(0,t)=u(1,t)=0.%
\end{array}%
\right.  \tag{5.2}
\end{equation}%
\bigskip

\bigskip 
\[
\begin{tabular}{|l|l|l|l|l|l|l|}
\hline
x & t & Exact Solution & Approximate Solution & Absolute Error & Relative
Error & Time \\ \hline
0.1 & 0.1 & $0.2938926262$ & $0.2938930965$ & $4.703\times 10^{-7}$ & $%
1.6002443\times 10^{-6}$ & $3.860$ \\ \hline
0.2 & 0.2 & $0.4755282582$ & $0.4755313577$ & $3.0995\times 10^{-6}$ & $%
6.518014327\times 10^{-6}$ & $3.016$ \\ \hline
0.3 & 0.3 & $0.4755282582$ & $0.4755183355$ & $9.9227\times 10^{-6}$ & $%
2.086668842\times 10^{-5}$ & $2.984$ \\ \hline
0.4 & 0.4 & $0.2938926261$ & $0.2939007109$ & $8.0848\times 10^{-6}$ & $%
2.750936663\times 10^{-5}$ & $3.000$ \\ \hline
0.5 & 0.5 & $0.0$ & $0.0000282140$ & $2.82140\times 10^{-5}$ & $\infty $ & $%
3.094$ \\ \hline
0.6 & 0.6 & $-0.2938926264$ & $-0.2939063137$ & $1.36873\times 10^{-5}$ & $%
4.657245120\times 10^{-5}$ & $3.031$ \\ \hline
0.7 & 0.7 & $-0.4755282583$ & $-0.4755305759$ & $2.3176\times 10^{-6}$ & $%
4.8737377\times 10^{-6}$ & $3.031$ \\ \hline
0.8 & 0.8 & $-0.4755282581$ & $-0.4755277748$ & $4.833\times 10^{-7}$ & $%
1.016343386\times 10^{-6}$ & $2.953$ \\ \hline
0.9 & 0.9 & $-0.2938926260$ & $-0.2938966580$ & $4.0320\times 10^{-6}$ & $%
1.371929624\times 10^{-5}$ & $3.204$ \\ \hline
1.0 & 1.0 & $0.0$ & $-3.690702068\times 10^{-7}$ & $3.690702068\times
10^{-7} $ & $\infty $ & $3.578$ \\ \hline
\end{tabular}%
\]

\begin{center}
\textbf{Table 1. }Numerical solutions for Example 5.1.

\[
\begin{tabular}{|l|l|l|l|l|l|}
\hline
$x$ & ES & AS & AE & RE & Time CPU(s) \\ \hline
$-0.80$ & $0.5877852522$ & $0.5877854278$ & $1.756\times 10^{-7}$ & $%
2.987485639\times 10^{-7}$ & $1.404$ \\ \hline
$-0.40$ & $0.9510565165$ & $0.9510565935$ & $7.70\times 10^{-8}$ & $%
8.096259125\times 10^{-8}$ & $1.373$ \\ \hline
$0.40$ & $-0.9510565165$ & $-0.9510565935$ & $7.70\times 10^{-8}$ & $%
8.096259125\times 10^{-8}$ & $1.373$ \\ \hline
$0.80$ & $-0.5877852522$ & $-0.5877854278$ & $1.756\times 10^{-7}$ & $%
2.987485639\times 10^{-7}$ & $1.404$ \\ \hline
\end{tabular}%
\]

\textbf{Table 2. }Numerical solutions for example 5.1 for t = 1.

\[
\begin{tabular}{|l|l|l|l|l|}
\hline
$x$ & AE [37] & AE [RKHSM] & RE [37] & RE [RKHSM] \\ \hline
$-0.80$ & $1.94E-05$ & $1.756E-07$ & $3.29E-05$ & $2.987485639E-7$ \\ \hline
$-0.40$ & $2.84E-07$ & $7.700E-08$ & $2.98E-07$ & $8.096259125E-8$ \\ \hline
$0.40$ & $2.84E-07$ & $7.700E-08$ & $2.98E-07$ & $8.096259125E-8$ \\ \hline
$0.80$ & $1.94E-05$ & $1.756E-07$ & $3.29E-05$ & $2.987485639E-7$ \\ \hline
\end{tabular}%
\]

\textbf{Table 3.} Comparison Absolute error and Relative error for Example
5.1.

\[
\begin{tabular}{|l|l|l|l|l|l|l|}
\hline
x & t & Exact Solution & Approximate Solution & Absolute Error & Relative
Error & Time \\ \hline
0.1 & 0.1 & $0.09917563307$ & $0.09917562$ & $1.307\times 10^{-8}$ & $%
1.317864035\times 10^{-7}$ & $1.688$ \\ \hline
0.2 & 0.2 & $0.1936096360$ & $0.1936060475$ & $3.5885\times 10^{-6}$ & $%
1.853471797\times 10^{-5}$ & $1.719$ \\ \hline
0.3 & 0.3 & $0.2794771925$ & $0.2794797547$ & $2.5622\times 10^{-6}$ & $%
9.167832184\times 10^{-6}$ & $0.593$ \\ \hline
0.4 & 0.4 & $0.3543825410$ & $0.354382187$ & $3.540\times 10^{-7}$ & $%
9.989205422\times 10^{-7}$ & $0.608$ \\ \hline
0.5 & 0.5 & $0.4173597173$ & $0.417357842$ & $1.8753\times 10^{-6}$ & $%
4.493246287\times 10^{-6}$ & $0.577$ \\ \hline
0.6 & 0.6 & $0.4685399032$ & $0.468538681$ & $1.2222\times 10^{-6}$ & $%
2.60852916\times 10^{-6}$ & $0.577$ \\ \hline
0.7 & 0.7 & $0.5087309794$ & $0.508730850$ & $1.294\times 10^{-7}$ & $%
2.543584040\times 10^{-7}$ & $0.593$ \\ \hline
0.8 & 0.8 & $0.5390654126$ & $0.539066210$ & $7.974\times 10^{-7}$ & $%
1.47922679\times 10^{-6}$ & $0.608$ \\ \hline
0.9 & 0.9 & $0.5607645949$ & $0.560770310$ & $5.7151\times 10^{-6}$ & $%
1.019162061\times 10^{-5}$ & $0.546$ \\ \hline
1.0 & 1.0 & $0.5750061826$ & $0.575005129$ & $1.0536\times 10^{-6}$ & $%
1.832328124\times 10^{-6}$ & $0.561$ \\ \hline
\end{tabular}%
\]%
\qquad \textbf{Table 4. }Numerical solutions for Example 5.2.

\[
\begin{tabular}{|l|l|l|l|l|l|}
\hline
$x$ & ES & AS & AE & RE & Time CPU(s) \\ \hline
$-0.80$ & $2.568109722$ & $2.568109726$ & $4\times 10^{-9}$ & $%
1.557565849\times 10^{-9}$ & $0.733$ \\ \hline
$-0.40$ & $2.985843344$ & $2.985843344$ & $0.0$ & $0.0$ & $0.686$ \\ \hline
$0.00$ & $3.141592654$ & $3.141592654$ & $0.0$ & $0.0$ & $0.702$ \\ \hline
$0.40$ & $2.985843344$ & $2.985843344$ & $0.0$ & $0.0$ & $0.686$ \\ \hline
$0.80$ & $2.568109722$ & $2.568109726$ & $4\times 10^{-9}$ & $%
1.557565849\times 10^{-9}$ & $0.733$ \\ \hline
\end{tabular}%
\]

\textbf{Table 5.} Numerical solutions for example 5.2 for t = 1.

\[
\begin{tabular}{|l|l|l|l|l|}
\hline
$x$ & AE [37] & AE [RKHSM] & RE [37] & RE [RKHSM] \\ \hline
$-0.80$ & $1.53E-08$ & $4E-09$ & $5.96E-09$ & $1.557565849E-09$ \\ \hline
$-0.40$ & $3.54E-10$ & $0.0$ & $1.18E-10$ & $0.0$ \\ \hline
$0.00$ & $1.62E-10$ & $0.0$ & $5.15E-11$ & $0.0$ \\ \hline
$0.40$ & $3.54E-10$ & $0.0$ & $1.18E-10$ & $0.0$ \\ \hline
$0.80$ & $1.53E-08$ & $4E-09$ & $5.96E-09$ & $1.557565849E-09$ \\ \hline
\end{tabular}%
\]

\textbf{Table 6. }Comparison Absolute error and Relative error for example
5.2.
\end{center}

\textbf{6. Conclusion}

In this study, linear and nonlinear SG Equations were solved by RKHSM. We
described the method and used it in some test examples in order to show its
applicability and validity in comparison with exact and other numerical
solutions. The obtained results show that this approach can solve the
problem effectively and need few computations. The results are satisfactory.
The results that we obtained were compared with the results that were
obtained by [37]. Numerical experiments on test examples show that our
proposed schemes are of high accuracy, and support the theoretical results.
According to these results, it is possible to apply RKHSM to linear and
nonlinear differential equations with initial and boundary conditions. It
has been shown that the obtained results are uniform convergent and the
operator that was used is a bounded linear operator. Thus the studies that
were compared with this method and this study show that RKHSM can be apply
to high dimensional partial differential equations, integral equations and
fractional differential equations without any transformation or
discretization and good results can be obtained.

\bigskip

\textbf{References}

[1] Russell J. Scott, Report of the committee of waves, Report of the 7th
Meeting of the British Associationfor the Advancement of Science, Liverpool,
1838, 417--496.

[2] Russell J. Scott, Report on waves, Report of the 14th Meeting of the
British Association for the Advancement of Science, John Murray, London,
1844, 311--390.

[3] D. J. Korteweg and G. de Vries, On the change of form of long-waves
advancing in a rectangular canal,and on a new type of long stationary waves,
Philos Mag 39 (1895), 422--443.

[4] N. J. Zabusky and M. D. Kruskal, Interaction of solitons in a
collisionless plasma and the recurrence of initial states, Phys Rev Lett 15
(1965), 240--243.

[5] M. J. Ablowitz and P. A. Clarkson, Solitons, nonlinear evolution
equations and inverse scattering, London Mathematical Society Lecture Note
Series 149, Cambridge University Press, Cambridge, 1991.

[6] A. D. Polyanin and V. F. Zaitsev, Handbook of Nonlinear Partial
Differential Equations, Chapman \& Hall/CRC, Boca Raton, Fla, USA, 2004.

[7] R. Rajaraman, Solitons and Instantons, North-Holland, Amsterdam, The
Netherlands, 1982.

[8] A. Mohebbi and M. Dehghan, \textquotedblleft High-order solution of
one-dimensional sine-Gordon equation using compact finite difference and
DIRKN methods,\textquotedblright\ Mathematical and Computer Modelling, vol.
51, no. 5-6, pp. 537--549, 2010.

[9] A. Scott, Nonlinear Science: Emergence and Dynamics of Coherent
Structure, vol. 8 of Oxford Texts in Applied and Engineering Mathematics,
Oxford University Press, Oxford, UK, 2nd edition, 2003.

[10] T. Dauxois and M. Peyrard, Physics of Solitons, Cambridge University
Press, 2006.

[11] M. Dehghan and A. Shokri, \textquotedblleft A numerical method for
one-dimensional nonlinear sine-Gordon equation using collocation and radial
basis functions,\textquotedblright\ Numerical Methods for Partial
Differential Equations, vol. 24, no. 2, pp. 687--698, 2008.

[12] M. Dehghan and D. Mirzaei, \textquotedblleft The boundary integral
equation approach for numerical solution of the one-dimensional sine-Gordon
equation,\textquotedblright\ Numerical Methods for Partial Differential
Equations, vol. 24, no. 6, pp. 1405--1415, 2008.

[13] A. G. Bratsos, \textquotedblleft A third order numerical scheme for the
two-dimensional sine-Gordon equation,\textquotedblright\ Mathematics and
Computers in Simulation, vol. 76, no. 4, pp. 271--282, 2007.

[14] A. G. Bratsos, \textquotedblleft A numerical method for the
one-dimensional sine-Gordon equation,\textquotedblright\ Numerical Methods
for Partial Differential Equations, vol. 24, no. 3, pp. 833--844, 2008.

[15] M. Dehghan and A. Shokri, \textquotedblleft A numerical method for
solution of the two-dimensional sine-Gordon equation using the radial basis
functions,\textquotedblright\ Mathematics and Computers in Simulation, vol.
79, no. 3, pp. 700--715, 2008.

[16] M. Lakestani and M. Dehghan, \textquotedblleft Collocation and finite
difference-collocation methods for the solution of nonlinear Klein-Gordon
equation,\textquotedblright\ Computer Physics Communications, vol. 181, no.
8, pp. 1392--1401, 2010.

[17] M. Dehghan and F. Fakhar-Izadi, \textquotedblleft The spectral
collocation method with three different bases for solving a nonlinear
partial differential equation arising in modeling of nonlinear
waves,\textquotedblright\ Mathematical and Computer Modelling, vol. 53, no.
9-10, pp. 1865--1877, 2011.

[18] A. H. A. Ali, \textquotedblleft Chebyshev collocation spectral method
for solving the RLW equation,\textquotedblright\ International Journal of
Nonlinear Science, vol. 7, no. 2, pp. 131--142, 2009.

[19] M.L. Ma, Z.M. Wu, A numerical method for one-dimensional nonlinear
Sine--Gordon equation using multiquadric quasi-interpolation, Chin. Phys. B
18 (2009) 3099--3103.

[20] A.Q.M. Khaliq, B. Abukhodair, Q. Sheng, M.S. Ismail, A
predictor-corrector scheme for the sine-Gordon equation, Numer. Methods
Partial Differ. Equations 16 (2000) 133--146.

[21] J. Chen, Z. Chen, S. Cheng, Multilevel augmentation methods for solving
the sine-Gordon equation, J. Math. Anal. Appl. 375 (2011) 706--724.

[22] N. Aronszajn, Theory of reproducing kernels, Trans. Amer. Math. Soc. 68
(1950) 337-404.

[23] M. Cui, Y. Lin, Nonlinear Numerical Analysis in the Reproducing Kernel
Spaces, Nova Science Publishers, New York, 2009.

[24] F.Geng, M. Cui, Solving a nonlinear system of second order boundary
value problems, J. Math. Anal. Appl. 327 (2007) 1167-1181.

[25] F.Geng, A new reproducing kernel Hilbert space method for solving
nonlinear fourth-order boundary value problems, Appl. Math. Comput. 213
(2009) 163-169.

[26] F.Geng, M. Cui, New method based on the HPM and RKHSM for solving
forced Duffing equations with integral boundary conditions, J. Comput. Appl.
Math. 233 (2009) 165-172.

[27] F.Geng, M. Cui, B. Zhang, Method for solving nonlinear initial value
problems by combining homotopy perturbation and reproducing kernel hilbert
spaces methods, Nonlinear Analysis: Real World Appl. 11 (2010) 637-644.

[28] F.Geng, M. Cui, Homotopy perturbation-reproducing kernel method for
nonlinear system of second order boundary value problems, J. Comput. Appl.
Math. 235 (2011) 2405-2411.

[29] F.Geng, M. Cui, A novel method for nonlinear two-point boundary value
problems: Combination of ADM and RKM, Appl. Math. Comput. 217 (2011)
4676-4681.

[30] M. Mohammadi, R. Mokhtari, Solving the generalized regularized long
wave equation on the basis of a reproducing kernel space, J. Comput. Appl.
Math. 235 (2011) 4003-4014.

[31] W. Jiang, Y. Lin, Representation of exact solution for the
time-fractional telegraph equation in the reproducing kernel space, Commun.
Nonlinear Sci. Numer. Simulat. 16 (2011) 3639-3645.

[32] Y. Wang, L. Su, X. Cao, X. Li, Using reproducing kernel for solving a
class of singularly perturbed problems, Comp. Math. Applic. 61 (2011)
421-430.

[33] B.Y. Wu, X.Y. Li, A new algorithm for a class of linear nonlocal
boundary value problems based on the reproducing kernel method, Appl. Math.
Letters 24 (2011) 156-159.

[34] H. Yao, Y. Lin, New algorithm for solving a nonlinear hyperbolic
telegraph equation with an integral condition, Int. J. Numer. Meth.
Biomedical Eng. 27 (2011) 1558-1568.

[35] F. Geng, M. Cui, A reproducing kernel method for solving nonlocal
fractional boundary value problems, Appl. Math. Comput. in press.

[36] Y. Wang, B. Wang, High-order multi-symplectic schemes for the nonlinear
Klein--Gordon equation, Appl. Math. Comput. 166 (2005) 608--632.

[37] Murat Sari and G%
\"{}%
urhan G%
\"{}%
urarslan. A sixth-order compact finite difference method for the
one-dimensional sine-Gordon equation. Int. J. Numer. Methods Biomed. Eng.,
27(7):1126--1138, 2011.

\end{document}